\documentclass[12pt,a4paper]{amsart2000}
\usepackage{amsmath2000,amssymb,amscd2000}

\newcommand{\ga}{\alpha}
\newcommand{\gl}{\lambda}
\newcommand{\ZZ}{\mathbb Z}
\newcommand{\wh}[1]{\Hat{#1}}
\newcommand{\cL}{\mathcal L}
\theoremstyle{plain}
\newtheorem{theorem}{Theorem}[section]
\newtheorem{lemma}[theorem]{Lemma}
\newtheorem{proposition}[theorem]{Proposition}

\theoremstyle{remark}

\begin{document}
\title[Moduli Spaces of Abelian Varieties]{An Isomorphism between
 Moduli Spaces of Abelian Varieties}

\author{Christina Birkenhake}
\author{Herbert Lange}
\address{Ch. Birkenhake\\Universit\"at Mainz\\Fachbereich Mathematik\\Staudingerweg
$9$\\D-$55099$ Mainz}
\address{H. Lange\\Mathematisches Institut\\
              Universit\"at Erlangen-N\"urnberg\\
              Bismarckstra\ss e $1\frac{ 1}{2}$\\
              D-$91054$ Erlangen\\
              Germany}
\email{birken@mathematik.uni-mainz.de}
\email{lange@mi.uni-erlangen.de}
\thanks{Supported by DFG Contracts Ba 423/8-1 and HU 337/5-1}
\maketitle






\begin{abstract}
In a previous paper we showed that for every polarization on an abelian
variety there is a dual polarization on the dual abelian variety. In this
note we extend this notion of duality to families of polarized abelian
varieties. As a main consequence we obtain an involution on the set
of moduli spaces of polarized abelian varieties of dimension $g$. In particular,
the moduli spaces $\mathcal
 A_{(d_1,\ldots,d_g)}$ and 
$\mathcal
 A_{(d_1,\frac{d_1d_g}{d_{g-1}},\ldots,\frac{d_1d_g}{d_2},d_g)}$ 
are canonically isomorphic.
\end{abstract}

\section{Introduction}

Let $k$ be an algebraically closed field and $(d_1,\ldots,d_g)$ a vector of
positive integers such that $d_i|d_{i+1}$ for all $i=1,\ldots,g-1$ and
$\mathrm{char}\,k\!\not\!|\,\,d_g$. Then the coarse moduli space
$\mathcal A_{(d_1,\ldots,d_g)}$ of polarized abelian varieties 
of type $(d_1,\ldots,d_g)$
over $k$
exists and is a quasi-projective variety of dimension $\frac{1}{2}g(g+1)$ over
$k$. The main result of this paper is the following theorem.
\begin{theorem}
There is a canonical isomorphism of coarse moduli spaces 
$$\mathcal
A_{(d_1,\ldots,d_g)}\stackrel{\sim}{\longrightarrow} \mathcal
A_{(d_1,\frac{d_1d_g}{d_{g-1}},\ldots,\frac{d_1d_g}{d_2},d_g)}.$$
\end{theorem}

For the proof we show that there is a canonical isomorphism of the corresponding
moduli functors. Namely, for every polarized projective 
abelian scheme $(A\longrightarrow
S,\lambda)$ we define a dual abelian scheme $(\wh{A}\longrightarrow
S,\lambda^\delta)$. Here $\wh{A}\longrightarrow S$ is the usual dual
abelian scheme and $\lambda^{\delta}$ is a polarization on $\wh{A}/S$ with
the property $(\lambda^{\delta})^{\delta}=\lambda$. Hence it makes sense to call
$\lambda^\delta$ the \textit{dual polarization} of $\wh{A}/S$.

There are several possibilities to define the polarization $\lambda^\delta$. In
section 2 we define a polarization $\lambda^D$ by inverting the isogeny
$\lambda:A\longrightarrow\wh{A}$. If $\lambda$ is given by a line bundle
$\mathcal L$ we define in section 4, using Mukai's Fourier functor, a relatively
ample line bundle $\wh{\mathcal L}$ on $\wh{A}$
inducing a polarization
$\phi_{\wh{\mathcal L}}$ of $\wh{A}$. Both polarizations are defined
even in the non separable case. We show that they are multiples of each other. 
If $\gl$ is a separable polarization of type 
$(d_1,\ldots,d_g)$, then $\lambda^{\delta}=d_1\gl^D$, and if moreover
$\gl=\phi_\cL$ with a relatively ample line bundle $\cL$ on $A/S$, then 
$\phi_{\wh{\mathcal L}}=d_1\cdots d_{g-1}\phi_\cL^D=d_2\cdots d_{g-1}\gl^{\delta}$.


\section{The Dual Polarization $\lambda^\delta$}

Let $\pi:A\longrightarrow S$ be an abelian scheme over a connected Noetherian
scheme $S$. Grothendieck showed in \cite{Groth}
that if $A/S$ is projective, the relative Picard functor $\mathcal{P}ic^0_{A/S}$
is represented by a projective abelian scheme
$\wh{\pi}:\wh{A}\longrightarrow S$, called the dual abelian scheme.
According to a theorem of Cartier and Nishi there is a canonical isomorphism
$A\stackrel{\sim}{\longrightarrow}\wh{\wh{A}}$ over $S$. We always
identify $A$ with its bidual $\wh{\wh{A}}$. \\
For any $a\in A(S)$ let $t_a:A\longrightarrow A$ denote the translation by $a$
over $S$. Any line bundle $\mathcal L$ on $A$ defines a homomorphism
$$
\phi_\cL:A\longrightarrow\wh{A}
$$
over $S$ characterized by $\phi_\cL(a):=t_a^*\cL\otimes\cL^{-1}$ for all $a\in
A(S)$. The biduality implies that $\phi_\cL$ is symmetric:
$\wh{\phi}_\cL=\phi_\cL$. Clearly $\phi_{\cL\otimes\pi^*\mathcal M}=\phi_\cL$
for any line bundle $\mathcal M$ on $S$. Moreover $\phi_\cL$ is an isogeny if
$\cL$ is relatively ample.\\

A \textit{polarization} of $\pi:A\longrightarrow S$ is a homomorphism
$\lambda:A\longrightarrow\wh{A}$ over $S$ such that for every geometric point
$s$ of $S$ the induced map $\lambda_s:A_s\longrightarrow\wh{A}_s$ is
of the form $\lambda_s=\phi_L$ for some ample line bundle $L$ of $A_s$. (Here
$\wh{A}_s:=(\wh{A})_s=(A_s)\wh{}$ )\\

Obviously for any relatively ample line bundle $\cL$ on $A/S$ the homomorphism
$\phi_\cL$ is a polarization. Conversely not every polarization of $A$ is of the
form $\phi_\cL$. However, if $\lambda:A\longrightarrow\wh{A}$ is a polarization,
then $2\lambda=\phi_\cL$ for some relatively ample line bundle $\cL$ on $A/S$
(see \cite{GIT} Prop 6.10, p.121).\\

Let $\lambda:A\longrightarrow\wh{A}$ be a polarization of $A$. Its kernel
$K(\lambda)$ is a commutative group scheme, finite and flat over $S$. According
to a Lemma of Deligne (see \cite{ot} p. 4) there is a positive integer $n$ such that
$K(\lambda)$ is contained in the kernel
$A_n:=\mathrm{ker}\{n_A:A\longrightarrow A\}$ where $n_A$ denotes
multiplication by $n$. The \textit{exponent} $e=e(\lambda)$ of the polarization
$\lambda$ is by definition the smallest such positive integer.

\begin{theorem}
\label{prop2.1}
There is a polarization $\lambda^D:\wh{A}\longrightarrow\wh{\wh{A}}=A$ of
$\wh{A}$, uniquely determined by $\lambda$, such that
$\lambda^D\circ\lambda=e_A$ and $\lambda\circ\lambda^D=e_{\wh{A}}$.
\end{theorem}
\begin{proof}    
Since $K(\lambda)\subset A_e$ by definition of the exponent, there is a
uniquely determined homomorphism $\lambda^D:\wh{A}\longrightarrow A$ such that
$\lambda^D\circ\lambda=e_A$. Hence
$$
(\lambda\circ\gl^D)\circ\gl=\gl\circ(\gl^D\circ\gl)=\gl\circ
e_A=e_{\wh{A}}\circ\gl.
$$
So $\gl\circ\gl^D=e_{\wh{A}}$, the homomorphism $\gl$
being surjective. It remains to show that $\gl^D$ is a polarization. For this
suppose $s$ is a geometric point of $S$. So $A_s$ and $\wh{A}_s$ are abelian
varieties over an algebraically closed field. It is well known that a
homomorphism $\varphi:\wh{A}_s\longrightarrow\wh{\wh{A}}_s=A_s$ is of the form
$\phi_M$ for some $M\in\mathrm{Pic}(\wh{A}_s)$ if and only if
$\varphi$ is symmetric: $\wh{\varphi}=\varphi$. So using
$\wh{\lambda_s}=\lambda_s$ and dualizing equation
$\gl^D_s\circ\gl_s=e_{A_s}$ give
$\gl_s\circ\wh{\gl^D_s}=\wh{\gl}_s\circ\wh{\gl^D_s}=e_{\wh{A}_s}$ Comparing
this with the equation $\gl_s\circ\gl^D_s=e_{\wh{A}_s}$ and using the fact that
$\gl_s$ and $e_{\wh{A}_s}$ are isogenies, implies $\wh{\gl^D_s}=\gl^D_s$. Hence
there is a line bundle $M$ on $\wh{A}_s$ such that $\gl^D_s=\phi_M$, and it
remains to show that $M$ is ample. For this note that
$\phi_{\gl_s^*M}=\wh{\gl}_s\circ\phi_M\circ\gl_s=\gl_s\circ\gl_s^D\circ\gl_s=e\gl_s.$
Now $e\gl_s$ being a polarization implies that $\gl_s^*M$, and thus also $M$ is
ample.
\end{proof}

A polarization $\lambda:A\longrightarrow \wh{A}$ is called 
\textit{separable} if for
every geometric point $s:\mathrm{Spec}(k)\longrightarrow S$ 
the characteristic of
$k$ does not divide the exponent $e(\lambda_s)$ of $\lambda_s$. In this case
$$
K(\lambda_s)=\mathrm{ker}(\lambda_s)\simeq(\ZZ/d_1\ZZ\times\cdots\ZZ/d_{g}\ZZ)^2
$$
with positive integers $d_1,\ldots,d_g$ such that $d_i|d_{i+1}$ and $d_g=e(\lambda_s)$
(see \cite{M} p. 294). Since the base scheme $S$ is connected and $K(\lambda)$
is an algebraic group scheme over $S$ with $K(\lambda)_s=K(\lambda_s)$ the integers
$d_1,\ldots,d_g$ do not depend on $s$. The vector $(d_1,\ldots,d_g)$ is called the
 \textit{type of the polarization} $\lambda$.

\begin{proposition}\label{prop2.2}
Suppose $\lambda$ is a separable polarization of the abelian scheme $A/S$ of
type $(d_1,\ldots,d_g)$. Then the polarization $\lambda^D$ of $\wh{A}/S$ is
separable of type $(1,\frac{d_g}{d_{g-1}},\ldots,\frac{d_g}{d_1})$.
\end{proposition}

\begin{proof}   
Let  $s$ be a geometric point of $S$. Equation
$\lambda^D_s\circ\lambda_s=e(\lambda_s)_{A_s}$ implies that the exponent
$e(\lambda_s^D)$ of $\lambda_s^D$ divides $e(\lambda_s)=d_g$, hence 
$\lambda^D$ is separable. Moreover
\begin{align*}
K(\lambda_s^D)=&\mathrm{ker}(\lambda^D_s)=\mathrm{ker}(e(\lambda_s)_{A_s})/\mathrm{ker}(\lambda_s)\\
\simeq&
(\ZZ/d_g\ZZ)^{2g}/(\ZZ/d_1\ZZ\times\cdots\times\ZZ/d_g\ZZ)^2\\
\simeq&\textstyle(\ZZ/\frac{d_g}{d_{g-1}}\ZZ\times\cdots\times\ZZ/\frac{d_g}{d_1}\ZZ)^2.
\end{align*}
\end{proof}

If $\lambda$ is a separable polarization of type $(d_1,\ldots,d_g)$, then Proposition 
\ref{prop2.2} implies that $(\lambda^D)^D$ is of type $(1,\frac{d_2}{d_1},\ldots,\frac{d_g}{d_1})$
and hence does not coincide with $\lambda$ if $d_1\not= 1$. However for
$$
\lambda^\delta:=d_1\lambda^D
$$
we have

\begin{proposition}\label{prop2.3}
$(\lambda^\delta)^\delta=\lambda$
\end{proposition}
\begin{proof}   
By definition $\gl^\delta$ is of type $(d_1,\frac{d_1d_g}{d_{g-1}},\ldots,\frac{d_1d_g}{d_2},d_g)$.
So both $\gl$ and $\gl^\delta$ have exponent $d_g$.
Applying Theorem \ref{prop2.1} to $\lambda$ and  $\lambda^\delta$ we get
$\lambda^\delta\circ\lambda=d_1(\gl^D\circ\gl)=(d_1d_g)_A=d_1(\gl^\delta\circ(\gl^\delta)^D)=
\lambda^\delta\circ(\lambda^\delta)^\delta$.
This implies the assertion.
\end{proof}

\section{Applications to Moduli Spaces}

Let $k$ be an algebraically closed field and $\mathcal S$ the category
of schemes of finite type over $k$. Fix a vector $(d_1,\ldots,d_g)$ of 
positive integers such that $d_i|d_{i+1}$ and $\mathrm{char}\,k\!\!\not|\,d_g$.
Consider the functor $\underline{\mathcal A}_{(d_1,\ldots,d_g)}:\mathcal S
\rightarrow \{sets\}$ defined by
$$
S\mapsto\{\text{isomorphism classes of pairs }(A\rightarrow S,\gl)\}
$$
where $A\rightarrow S$ is a projective abelian scheme of relative dimension
$g$ over $S$, and $\gl:A\rightarrow\wh{A}$ is a polarization of type 
$(d_1,\ldots,d_g)$. 
Note that any such polarization $\gl$ of $A$ is separable by the assumption 
on the characteristic $k$ of $S$.

Recall that a \textit{coarse moduli scheme}
for the functor $\underline{\mathcal A}_{(d_1,\ldots,d_g)}$ is a scheme
$\mathcal A_{(d_1,\ldots,d_g)}$ of finite type over $k$ 
admitting a morphism of functors 
$$
\ga:\underline{\mathcal A}_{(d_1,\ldots,d_g)}
(\,\cdot\,)\rightarrow\mathrm{Hom}(\,\cdot\,,\mathcal A_{(d_1,\ldots,d_g)})
$$
such that
\begin{enumerate}
\item $\ga(\mathrm{Spec}\,k):\underline{\mathcal A}_{(d_1,\ldots,d_g)}(\mathrm{Spec}\,k)
\rightarrow
\mathbf{Hom}(\mathrm{Spec}\,k,\mathcal A_{(d_1,\ldots,d_g)})$ is a bijection, and
\item for any morphism of functors $\beta:\underline{\mathcal A}_{(d_1,\ldots,d_g)}
(\,\cdot\,)\rightarrow\mathrm{Hom}(\,\cdot\,, B)$ there is
a morphism of schemes $f:\mathcal A_{(d_1,\ldots,d_g)}\rightarrow  B$
such that $\beta=\mathrm{Hom}(\,\cdot\,,f)\circ\ga$.
\end{enumerate}

It is well known that a uniquely determined coarse moduli scheme 
$\mathcal A_{(d_1,\ldots,d_g)}$ exists (see \cite{GIT}). Moreover it is clear 
from the definition that an isomorphism of functors induces an isomorphism of 
the corresponding coarse moduli schemes. This 
will be used to prove the following 

\begin{theorem}\label{thm4.1}
There is a canonical isomorphism of coarse moduli schemes
$$
\mathcal A_{(d_1,\ldots,d_g)}\longrightarrow
\mathcal A_{(d_1,\frac{d_1d_g}{d_{g-1}},\ldots,\frac{d_1d_g}{d_2},d_g)}.
$$
\end{theorem}
\begin{proof}   
By what we have said above it suffices 
to show that there is a canonical isomorphism of functors
$\underline{\mathcal A}_{(d_1,\ldots,d_g)}\rightarrow
\underline{\mathcal A}_{(d_1,\frac{d_1d_g}{d_{g-1}},
\ldots,\frac{d_1d_g}{d_2},d_g)}$. But
for any $S\in\mathcal S$ the canonical map 
$\underline{\mathcal A}_{(d_1,\ldots,d_g)}(S)\rightarrow
\underline{\mathcal A}_{(d_1,\frac{d_1d_g}{d_{g-1}},
\ldots,\frac{d_1d_g}{d_2},d_g)}(S)$, 
defined by 
$$
(A\rightarrow S,\gl)\mapsto(\wh{A}\rightarrow S,\gl^\delta)
$$
has an inverse, since $\wh{\wh{A}}=A$ by the biduality theorem and
$(\gl^\delta)^\delta=\gl$ by Proposition \ref{prop2.3}. Moreover 
this is an isomorphism of functors, since $(\gl_T)^\delta=(\gl^\delta)_T$
for any morphism $T\rightarrow S$ in $\mathcal S$.
\end{proof}

\section{The Dual Polarization via The Fourier Transform}
In Section 2 we defined for every polarization $\lambda$ of a projective
 abelian
scheme $A/S$ a polarization $\lambda^D$ of the dual abelian scheme $\wh{A}/S$.
Now suppose $\lambda=\phi_\cL$ for a relatively ample line bundle $\cL$ on
$A/S$. In this section we apply Mukai's Fourier transform to define a
relatively ample line bundle $\wh{\cL}$ on $\wh{A}/S$ which induces a multiple
of $\lambda^D$.\\

Let $\pi:A\longrightarrow S$ be a projective abelian scheme over a connected
Noetherian scheme $S$ of relative dimension $g$. Let $\mathcal P=\mathcal P_A$
denote the normalized Poincar\'e bundle on $A\times_S\wh{A}/S$. Here normalized
means that both $(\epsilon\times1_{\wh{A}})^*\mathcal P$ and
$(1_A\times\wh{\epsilon})^*\mathcal P$ are trivial, where
$\epsilon:S\longrightarrow A$ and $\wh{\epsilon}:S\longrightarrow \wh{A}$ are
the zero sections. Denote by $p_1$ and $p_2$ the projections of
$A\times_S\wh{A}$. A coherent sheaf $\mathcal M$ on $A$ is called
\textrm{WIT}-\textit{sheaf of index} $i$ on $A$ if $R^j{p_2}_*(\mathcal
P\otimes p_1^*\mathcal M)=0$ for $j\not =i$. In this case
$$
F(\mathcal M):=R^i{p_2}_*(\mathcal P\otimes p_1^*\mathcal M)
$$
is called the \textit{Fourier transform} of $\mathcal M$. Let $\cL$ be a
relatively ample line bundle on $A/S$. The scheme $S$ being connected, the
degree $d:=h^0(\cL_s)$ is constant for all geometric points $s$ of $S$, it is
called the \textit{degree} of $\cL$. By the base change theorem $\cL$ is a
\textrm{WIT}-sheaf of index $0$ and its Fourier transform $F(\cL)$ is a vector
bundle of rank $d$. Define
$$
\wh{\cL}:=\Bigl(\det{F(\cL)}\Bigr)^{-1}.
$$
Let $e=e(\phi_\cL)$ be the exponent of the polarization $\phi_\cL$. The main
result of this section is the following
\begin{theorem}\label{thm3.1}
$\wh{\cL}$ is a relatively ample line bundle on $\wh{A}/S$ such that
$$
e\phi_{\wh{\cL}}=d\phi_\cL^D.
$$
\end{theorem}
If $\phi_\cL$ is a separable polarization of type $(d_1,\ldots,d_g)$ then
$e=d_g$ and $d=d_1\cdot\ldots\cdot d_g$. So $\phi_{\wh{\cL}}$ is a multiple of
$\phi_\cL^D$:
$$
\phi_{\wh{\cL}}=d_1\cdot\ldots\cdot d_{g-1}\phi_\cL^D.
$$
In the non separable case $e|d^2$ by Deligne's Lemma but it is not clear to us
whether $e|d$.

For the proof we need some preliminaries.
\begin{lemma}\label{lem3.2}
Let $\varphi:A\longrightarrow B$ be an isogeny of abelian schemes and $\mathcal
M$ a \textrm{WIT}-sheaf on $A$, then $\varphi_*\mathcal M$ is a
\textrm{WIT}-sheaf of the same index on $B$ and
$$
F_B(\varphi_*\mathcal M)=\wh{\varphi}\,^*F_A(\mathcal M).
$$
Here $F_A$ and $F_B$ denote the Fourier transforms of $A$ and $B$ respectively.
\end{lemma}
The proof is the same as in the absolute case (see \cite{M1} (3.4)).

\begin{lemma}\label{lem3.3}
Let $\mathcal M$ be a \textrm{WIT}-sheaf of index $i$ on $A$ and $\mathcal F$ a
locally free sheaf on $S$. Then $\pi^*\mathcal F\otimes \mathcal M$ is a
\textrm{WIT}-sheaf of index $i$ on $A$ and
$$
F(\pi^*\mathcal F\otimes \mathcal M)=\wh{\pi}^*\mathcal F\otimes F(\mathcal M).
$$
\end{lemma}
\begin{proof}   
Using $\pi\circ
p_1=\wh{\pi}\circ p_2$ and the projection formula we get
\begin{align*}
R^j{p_2}_*(\mathcal P\otimes p_1^*\pi^*\mathcal F\otimes p_1^*\mathcal M)
   =&   R^j{p_2}_*(p_2^*\wh{\pi}^*\mathcal F\otimes\mathcal P\otimes p_1^*\mathcal M)\\
   =&\wh{\pi}^*\mathcal F\otimes  R^j{p_2}_*(\mathcal P\otimes p_1^*\mathcal M)\\
  =&\begin{cases}
                 \wh{\pi}^*\mathcal F\otimes F(\mathcal M) & \text{if $j=i$}\\
                 0                                         & \text{otherwise}.
    \end{cases}
\end{align*}
This implies the assertion.
\end{proof}

\begin{proposition}\label{prop3.4}
Let $\cL$ be a relatively ample line bundle on $A/S$.
Then $\cL^{-1}$ is a \textrm{WIT}-sheaf of index $g$ satisfying
$\phi_\cL^*F(\cL)=\pi^*\pi_*\cL\otimes\cL^{-1}$.
\end{proposition}

\begin{proof}   
The first assertion follows from Serre duality.
Denote by $q_1,q_2:A\times_SA\rightarrow A$ the projections.
Then
\begin{align*}
\pi^*\pi_*\cL\otimes\cL^{-1}&={q_2}_*(q_2-q_1)^*\cL\otimes\cL^{-1}\\
 &\hspace{1,2cm}\Bigl(\text{(using flat base change with}\quad\pi\circ
q_2=\pi\circ(q_2-q_1)\Bigr)\\
 &={q_2}_*\Bigl((q_2-q_1)^*\cL\otimes q_2^*\cL^{-1}\Bigr)\\
 &={q_2}_*\Bigl(((-1)_A\times\phi_\cL)^*\mathcal P\otimes q_1^*(-1)_A^*\cL\Bigr)\\
 \intertext{\hspace{1cm}$ \Bigl($applying the formula $((-1)_A\times\phi_\cL)^*\mathcal P=(q_2-q_1)^*\cL\otimes
 q_1^*(-1)_A^*\cL^{-1}\otimes q_2^*\cL^{-1}\Bigr)$}
 &= {q_2 }_*\Bigl((-1)_A\times\phi_\cL)^*(\mathcal P\otimes p_1^*\cL)\Bigr)\\
 &=\phi_\cL^*{p_2 }_*(\mathcal P\otimes p_1^*\cL)\quad\text{  (by flat base change)}\\
 &=\phi_\cL^*F(\cL).
 \end{align*}
\end{proof}

\begin{proof}[of Theorem \ref{thm3.1}]
By Proposition \ref{prop3.4} we have
$$
\phi_\cL^*(\wh{\cL})=\Bigl(\det\phi_\cL^*F(\cL)\Bigr)^{-1}=
\Bigl(\det(\pi^*\pi_*\cL\otimes\cL^{-1})\Bigr)^{-1}=\pi^*(\det\pi_*\cL)^{-1}\otimes\cL^d,
$$
since $\pi^*\pi_*\cL$ is a vector bundle of rank $d$.
This implies that $\wh{\cL}$ is relatively ample. Moreover
$$
\phi_\cL\circ\phi_{\wh{\cL}}\circ\phi_\cL=\phi_{\phi_\cL^*\wh{\cL}}
=\phi_{\pi^*(\det\pi_*\cL)^{-1}\otimes\cL^d}=\phi_{\cL^d}=d\phi_\cL
$$
implies that
$$
\phi_{\wh{\cL}}\circ\phi_\cL=d_A\quad\text{and}\quad\phi_\cL\circ\phi_{\wh{\cL}}=d_{\wh{A}}.
$$
Comparing this with Theorem \ref{prop2.1} gives the assertion.
\end{proof}

\end{document}